\documentclass[11pt,a4paper]{amsart}

\usepackage[latin1]{inputenc}
\usepackage[english]{babel}
\usepackage{amsmath}

\usepackage{amssymb, mathabx}
\usepackage[numbers]{natbib}
\usepackage{graphicx}
\usepackage{braket}
\usepackage[pdftex,plainpages=false,colorlinks,hyperindex,bookmarksopen,linkcolor=red,citecolor=blue,urlcolor=blue]{hyperref}

\usepackage{epstopdf}

\usepackage{braket}

\bibpunct{[}{]}{;}{n}{,}{,}
\usepackage[hmargin=3cm,vmargin={3.5cm,4cm}]{geometry}

\theoremstyle{definition}                                 %stile corsivo
\theoremstyle{definition}                           %stile roman
\theoremstyle{remark}                             %stile per osservazioni
\usepackage{color}

\def\d{\partial}

\newcommand{\ee}{\end{eqnarray}}
\def\d{\partial}

\def\pni{\par\noindent}
\def\vsh{\smallskip}

\def\vsp{\vsh\pni} %% ie. \smallskip + \par

 \def\G{{\mathcal{G}}}

\def\d{\partial}
\def\RR{\vbox {\hbox to 8.9pt {I\hskip-2.1pt R\hfil}}\;}

\numberwithin{equation}{section}

\allowdisplaybreaks

\begin{document}
\title{On the evolution of 	\\ fractional diffusive waves}

    \author{Armando Consiglio$^1$}
		\address{${}^1$ 
        Institut f\"{u}r Theoretische Physik und Astrophysik, 
        Universit\"{a}t W\"{u}rzburg, D-97074 W\"{u}rzburg, GERMANY.}
         \email{armconsiglio@gmail.com}  		
		
    \author{Francesco Mainardi$^2$}
    	    \address{${}^2$ Department of Physics $\&$ Astronomy, University of 	
    	    Bologna and INFN. Via Irnerio 46, I-40126 Bologna, ITALY.}
			\email{francesco.mainardi@bo.infn.it}

   % \date{December 2019}
    \maketitle 

\date{December 2019} 
\begin{abstract}
In physics,  phenomena of diffusion and wave propagation have great relevance; these physical processes are governed in the simplest cases  by partial differential equations of order 1 and 2 in time, respectively.
It is known that whereas
the diffusion equation describes a process where the disturbance spreads
infinitely fast, the propagation velocity of the disturbance is a constant
for the wave equation. 
By replacing the time  derivatives in the above standard equations with pseudo-differential operators interpreted as derivatives of non integer order (nowadays misnamed as of fractional order)  
% of order $\alpha$ in time with $1 \leq \alpha \leq 2$,
 we are lead to generalized processes of diffusion  that may be interpreted as slow diffusion and interpolating  between diffusion and wave propagation. In mathematical physics, we may refer these interpolating processes  to as fractional diffusion-wave phenomena. The use of the Laplace transform in the analysis of the Cauchy and Signalling problems leads to  special functions of the Wright type.
 %% nowadays known as M-Wright function
 \\ 
In this work 
%we  show that the time-fractional diffusion-wave equation interpolates between the two different %responses, 
we analyze  and simulate both the situations in which the input  function is a Dirac delta generalized function and a box function, restricting ourselves to the Cauchy problem.
 In the first case we get  the fundamental solutions
(or Green functions) of the problem whereas  in the latter case the solutions are obtained 
by a space convolution of the Green function with the input  function.
  In order to clarify the matter for the non-specialist readers, we briefly recall the basic and essential notions of  the fractional calculus (the mathematical theory that regards the integration and differentiation of non-integer order) with a look at the history of this discipline.  
%\keywords{Fractional calculus \and Diffusion-waves phenomena \and  Wright % functions \and Cauchy and Signaling Problems.}  
\end{abstract}
\vsp
\vsp
\vsp
{\bf Published on line (6 Dec. 2019) in Ricerche di Matematica, Springer
(without the photos in Fig. 1 and in Figs. 3-4); 
 DOI 10.1007/s11587-019-00476-6}
  %\end{abstract}
%%%%
\section{Introduction}
It is known that mathematical models are usually governed by differential  and/or integral  equations of integer order. However, we can extend 
integration and differentiation to  any order by using the tools of the so-called \emph{Fractional Calculus}, that in recent years has gained considerable interest also because of its applications in different fields including 
 Physics,  Engineering, Biology,
	Economics and Finance, Geophysics,
	 Computer Science
and so on...
\\
Fundamental phenomena of Physics, such as diffusion and wave propagation, are governed in their simplest form by the following partial differential equations  (PDE) that in a standard notation  read
%\medskip
\begin{itemize}
\item Diffusion equation: 
	\begin{equation}
	\frac{\partial u(\mathbf{r}, t)}{\partial t } = D \nabla ^2 u(\mathbf{r}, t),
	\end{equation}
	\item Wave equation:
	\begin{equation}
	\frac{\partial ^2 u(\mathbf{r}, t)}{\partial t ^2} = c^2 \nabla ^2 u(\mathbf{r}, t).
	\end{equation}
\end{itemize}
\vsp
The above PDE's can be extended by replacing   the time and/or  spatial derivatives of integer order  with some pseudo-differential operators interpreted as derivative of non integer order. In this paper we consider the generalizations of the above PDE's by replacing the time derivatives with derivatives of order $\alpha \in (0,2]$
and we restrict our attention  to a single space dimension $x$   
\vsp
     This work is organized as follows. 
 In Section 2, we show the essential mathematical notions  
 necessary to understand the integration and differentiation of non-integer order, that nowadays 
 are ascribed to the mathematical discipline referred generally to as {\it Fractional Calculus}. 
 In particular we point out   the approach of Professor Caputo  on which  Fractional Calculus has grown and found several applications. since the late 1960's.
\vsp
 Then, in Sections 3 and 4 we provide some historical notes distinguishing two different eras  with the advent of Caputo's approach. Indeed,  in Section 3 we recall
  the major mathematicians  who in the last centuries have contributed in the birth of development of the 
  Fractional Calculus whereas in Section 4 we outline how the Caputo approach  was popularized.
  \vsp
  The core of the paper is found in Section 5, where we  consider from a mathematical view point the evolution of the so called fractional diffusive waves generated by  partial differential equations 
  with time derivatives of non integer order.  
  \vsp
 Finally in Section 6, we  provide some concluding remarks paying attention to work to be done in the next future.
 
\section{Essentials of Fractional Calculus: Fractional Derivatives and Integrals}
The starting point for the Fractional Calculus is given by the \emph{Cauchy formula for repeated integration} where $a \in  \RR$.  Taking as  independent variable  the time $t$ and denoting  by $f(t)$
a function assumed to be sufficiently well-behaved for the next considerations, we recall the formula
for the $n$ repeated integral
\begin{equation}
{}_{a}I^{n}_{t} f(t):=  \int_{a}^t \int_{a}^{\tau_{n-1}} \!\cdots\!  \int_{a}^{\tau_{1}}
\! f(\tau)\, d\tau\, d\tau_1 \! \cdots \!  d\tau_{n-1} \!=\! \frac{1}{(n-1)!}
\int_{a}^{t} (t-\tau)^{n - 1}\,f(\tau)\, d\tau.
\end{equation}
The resulting function turns out to be the primitive of order $n$  of $f(t)$, that is vanishing in 
$t = a$ along with its derivatives of order $1, \cdots ,  n-1$.
The basic idea of the fractional calculus is in replacing  the factorial  term $(n-1)!$ with the corresponding representation by the Gamma function and then replacing $n$ with a positive real number $\alpha$. 
 So, writing $(n-1)! = \mathit{\Gamma} (n)$ in Eq. (3), it follows the definition of \emph{fractional integral} ${}_{a}I^{\alpha}_{t} f(t)$ of order $\alpha$  ($\alpha > 0$)
\begin{equation}
{}_{a}I^{\alpha}_{t} f(t): = 
\frac{1}{\mathit{\Gamma} (\alpha)}\int_{a}^{t}{(t-\tau)^{\alpha - 1}\,f(\tau)\, d\tau},
\end{equation}
where now $\alpha$ is no longer restricted only to positive, integer values.
\vsp
The \emph{fractional derivative} of order $\alpha$ is so defined (in the Riemann-Liouville sense) as the left inverse operator of the corresponding fractional integral and reads:
\begin{equation}
{}_{a}D^{\alpha}_{t} f(t) := \frac{d^n}{dt^n}{{}_{a}I^{n-\alpha}_{t} f(t)}, \quad  n = [\mathfrak (\alpha )] +1 .
\end{equation}
%\vsp
In view of our applications we take $ a=0$ so $t  \ge 0$ and, from now on, we agree to delete $a$ in the notations for fractional integrals and derivatives.
Henceforth, we  consider the time fractional derivative in the Caputo sense, 
defined  for $n-1<\alpha \le n$ as:
\begin{equation}
D^{\alpha}_{C} f(t): = I^{n-\alpha}_{t}\frac{d^n}{dt^n} f(t) =
\frac{1}{\mathit{\Gamma} (n-\alpha)}\int_{0}^{t}{\frac{f^{(n)}(\tau)}
{(t-\tau)^{\alpha + 1 -n}}\, d\tau}.
\end{equation}
This definition provides us a regularization of the Riemann-Liouville fractional derivative at $t=0$ 
as shown in 1997  by Gorenflo and Mainardi  \cite{Gorenflo-Mainardi CISM97}. 
\vsp
As a consequence of its definition,   the Caputo derivative is vanishing when
the derivative of integer order ${f^{(n)}(t)}$ is zero;  in particular 
\begin{equation}
D^{\alpha}_C 1 = 0, ~~~~~ \alpha > 0\,.
\end{equation}
The Caputo derivative 
 appears suitable to be treated by the Laplace Transform technique for causal systems, that is quiescent for $t<0$. In a standard notation we get  
\begin{equation}
D^{\alpha}_{C}f(t) \div s^{\alpha} \tilde{f}(s) - \sum_{k=0}^{n-1}{f^{(k)}(0^+)\, s^{\alpha - 1 - k}},~~ n-1 < \alpha \le n\,.
\end{equation}
It is requested the knowledge of the initial values of the function and
of its integer derivatives of order $k = 1, 2, ..., n-1$.as it is for  $\alpha = n$.
%%5    \newpage
\section{Historical notes before Caputo's era}
Fractional Calculus may be considered born on September 30$^{th}$, 1695 when 
L'Hopital asked to Leibniz what would be if $n= 1/2$ in his notation for the n$^{th}$ derivative, $\frac{d^n}{dx^n}$. %\\
Leibniz answered: ``An apparent paradox, from which one day useful consequences will be drawn."
\begin{figure}[h!]%
	\centering
	%\subfloat[\tiny{Gottfried Wilhelm (von) Leibniz}]{
	{{\includegraphics[width=2.5cm]{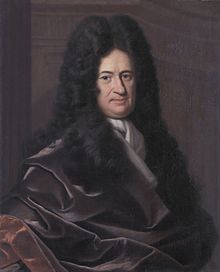} }}%
	\qquad
	%\subfloat[\tiny{Guillaume Fran\c{c}ois Antoine, Marquis de l'H\^{o}pital}]
	{{\includegraphics[width=2.5cm]{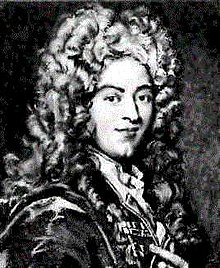} }}%
	\caption{LEFT: Gottfried Wilhelm (von) Leibniz;   RIGHT: Guillaume Fran\c{c}ois Antoine, Marquis de l'H\^{o}pital.}
	\label{fig:Leiniz-Hopital}%
\end{figure}
\vsp
From that time we get  a very long story involving eminent mathematicians who are pointed out in the next figure. More recently the term {fractional calculus} was used but it is a misnomer kept only for historical reasons related to the  discussion between
Marquis de l'H\^{o}pital and Leibniz.  
\vsp 
For more details see the  historical notes by the late  
Bertram Ross, the  organizer of the first conference devoted to fractional calculus in 1974 \cite{Ross History-BOOK75}, and, more recently, by Machado and Kiryakova \cite{Machado-Kiryakova_HFCA19}
from which the following picture is taken.
\begin{figure}[h!]%
	\centering
	\includegraphics[width=6.8cm]{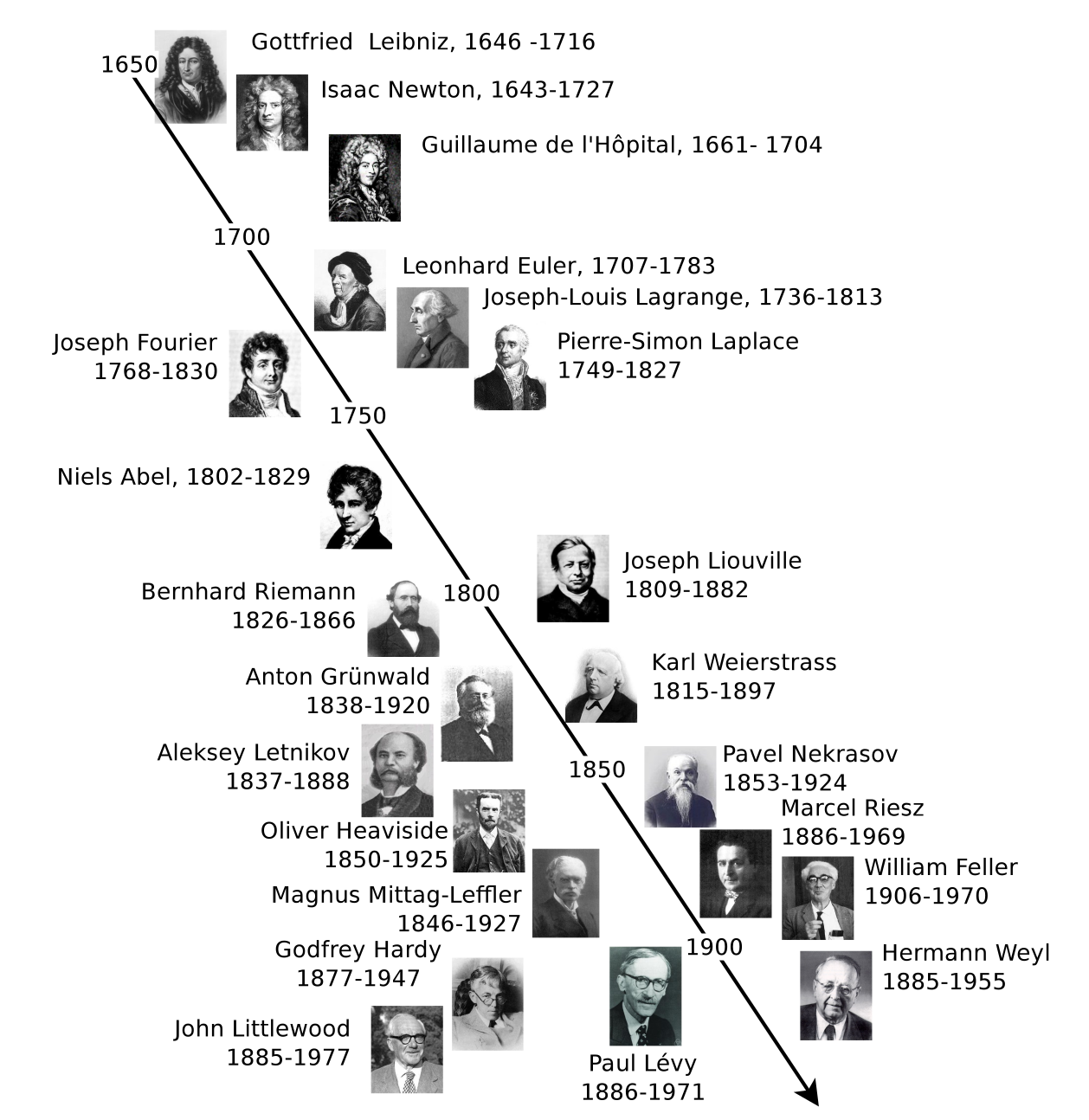}
	\caption{The major mathematicians in Fractional Calculus up to 1950's.} 
	\label{fig:story}%
\end{figure}
     
\section{Historical notes in  Caputo's era}

 In this paper, attention  is indeed devoted 
to the form (6) of fractional derivative alternative to that originally used by Liouville and Riemann expressed by Eq (5).
\vsp
The form (6), where  the orders of fractional integration and ordinary differentiation are interchanged,
is  nowadays known as the {\it Caputo derivative}.
   As a matter of fact, such a form is found in a paper  by Liouville himself
   as  noted
by Butzer and Westphal \cite{Butzer-Westphal 00} in the 2000  book edited by  R. Hilfer,
but Liouville, not recognizing its role,   disregarded this notion.
As far as we know, 
up to  the middle of the twentieth century, most authors did not take notice
of the difference between the two forms and of the possible
use of  the alternative form. 
  Even in the classical book on Differential and Integral Calculus published in English in 1936 by the eminent mathematician 
  R. Courant,
the two forms of the fractional derivative were considered as equivalent,
see   \cite{Courant BOOK36}, pp.  339-341.
\vsp 
In the  late sixties of the past century    the relevance of the alternative form 
was finally recognized. In fact, in  1968 the Soviet Scientists Dzherbashyan and Nersesyan 
  \cite{Dzherbashyan-Nersesyan 68} and then in 1989-90 Kochubei \cite{Kochubei 89,Kochubei 90}
  used the alternative form (6)  in dealing  with Cauchy problems
  for differential equations of fractional order.
  However, formerly since 1967,  Caputo \cite{Caputo 67,Caputo BOOK69}
  had introduced  this form as given by Eq. (6) proving the corresponding rule in the Laplace transform domain, as Eq. (8).  
   With his derivative   Caputo  was  thus able to  generalize the  rule for the Laplace transform  
 of a derivative of integer order and to solve some problems
 in Seismology in a proper way related to the dissipation function.
 Indeed.  it was the corresponding rule in the Laplace transform the basic and original  idea of the Caputo form of fractional derivative, not present in any previous paper. 
 Soon later, this derivative was adopted by Caputo and Mainardi %%(1971a),  (1971b), see
 \cite{Caputo-Mainardi 71PAG}, \cite{Caputo-Mainardi 71RNC}  
 in the framework of the theory of {\it Linear Viscoelasticity}.
 %% \cite{Caputo-Mainardi 71PAG,Caputo-Mainardi 71RNC}.
 %%%%%%%
  %%%%%%%%%%%%%%% THE END of FOOTNOTE
 \vsp
  Since the seventies of the past century a number of authors
  %% , very often  ignoring the works by   Caputo, Mainardi and Kochubei
  have re-discovered and used the alternative form, recognizing its major utility
  for solving physical problems with standard initial conditions.
  Although  several  papers by different authors
  appeared %% including Caputo and  Kochubei,
  where the alternative derivative
  was adopted, it was only in the late nineties,    with  the  tutorial paper   by Gorenflo and Mainardi
  in 1997 \cite{Gorenflo-Mainardi CISM97} and   the 1999  book by Podlubny 
   \cite{Podlubny BOOK99}, that such form was  popularized. 
  % However the first article where the Caputo derivative was present in a title was in 1998 by Luchko and %Gorenflo \cite{LuGo 99}.
%%%
  % In the above  references   the Caputo form  was  named the 
  Nowadays the term  {\it Caputo fractional derivative}
   is universally accepted in the literature. 
  The reader, however, is alerted that in a very few papers the Caputo derivative is referred 
to as the {Caputo--Dzherbashyan}  derivative. 
Note also the transliteration as Djrbashyan.
\vsp
  It should be noted that it was Mainardi, being initially  an Italian  researcher in Geophysics as a PhD student of Prof. Caputo, who pointed out the Caputo form (6) published  in a Geophysical journal \cite{Caputo 67}
  and in an Italian book \cite{Caputo BOOK69} to the attention of eminent colleagues in Mathematics 
  including 
  the late Rudolf Gorenflo 
  \cite{Gorenflo-Mainardi CISM97}, \cite{Gorenflo-Kilbas-Mainardi-Rogosin BOOK14},
  Yuri Luchko \cite{LuGo 99},  Igor Podlubny \cite{Podlubny BOOK99}, the late Anatoly Kilbas \cite{Kilbas-Srivastava-Trujillo BOOK06}. Kai Diethelm \cite{Diethelm  BOOK10}.
 Also Virginia Kiryakova, the Editor-in-Chief 
  of  the major journal devoted to Fractional Calculus, that is Fractional Calculus and Applied Analysis
    (FCAA), was alerted. 
  Indeed, the papers by Caputo in 1967 \cite{Caputo 67} and by Caputo and Mainardi in 1971
  \cite{Caputo-Mainardi  71PAG} were reprinted 
in FCAA in 2007 on the occasion of  the eighty birthday of Professor Caputo.
Just in this volume  there were the  foreword by Mainardi 
\cite{Mainardi FCAA07} 
to the special issue of FCAA dedicated to the 40 years of Caputo derivative and to the 80th anniversary of Prof Caputo. 
This issue contains  a tutorial survey by Mainardi and Gorenflo \cite{Mainardi-Gorenflo FCAA07}
and  the biographical data of Prof. Caputo and his photo reported  below.
\begin{figure}[h!]%
	\centering
	{{\includegraphics[width=3.5cm]{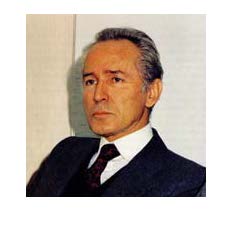} }}
	\caption{{Professor Michele Caputo in 1967.}}
	\label{fig:caputo}%
\end{figure}
%%%
%\\
\begin{figure}[h!]%
	\centering
	{{\includegraphics[width=3.5cm]{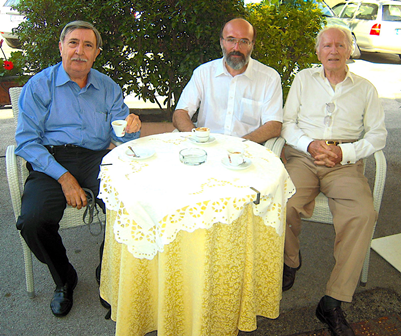} }}%
	\caption{{Mainardi, Podlubny and Caputo, Ravello, September 2012.}}
	\label{fig:mainardi-podlubny-caputo}%
\end{figure}
We enclose also the photo depicting Mainardi, Podlubny and Caputo at Ravello, where Mainardi and Podlubny delivered each a  one-week course on fractional calculus  and Caputo a seminar in the framework of the 37-th summer school on Mathematical Physics organized by the National Group of Mathematical Physics, held just in Ravello (17-29 September 2012) directed by Professors  Rionero and Ruggeri.
  %  We note that the {\it Caputo fractional  derivative} coincides with that
 % introduced,  independently and  a few later,
 % by  {\it Dzherbashyan \& Nersesyan} \cite{Dzherbashyan-Nersesyan 68}
% The role of this fractional derivative in treating initial value problems
 %  has been   recognized  also by  Kochubei,
 % see \eg \cite{Kochubei 89,Kochubei 90}, and later rediscovered  
%% This name was essentially introduced  by Podlubny who in his book
%% in order to distinguish it from the classical notion
%% usually refereed to Riemann and Liouville.  
%% We remind that the Laplace transform rule (1.30)
%% was practically the starting point of Caputo himself in defining
%% his generalized derivative in the late sixties, \cite{Caputo 67,Caputo BOOK69}.
%% Later, Caputo and Mainardi in 1971 \cite{Caputo-Mainardi 71PAG,Caputo-Mainardi 71RNC}
%% and Mainardi in the nineties, see e.g. \cite{Mainardi WASCOM93,Mainardi CHAOS96},
%% have followed the notation involving  a convolution with the so-called
%% Gel'fand-Shilov (generalized) function 
%% $\, \Phi_{\lambda}(t):= t_+^{\lambda -1}/\Gamma(\lambda)\,$
%% introduced  in \cite{Gelfand-Shilov BOOK64},
%% as  discussed at the end of Section 1.2, see Eqs. (1.31)-(1.35).
%% 

%\hskip-0.5truecm
%\newpage

\section{Fractional calculus in diffusion-wave problems}
The time fractional diffusion-wave equation reads: 
\begin{equation}
\frac{\partial ^{\alpha} u(\mathbf{r}, t)}{\partial t^{\alpha}} = D\nabla ^2 u(\mathbf{r}, t), ~~~ D > 0, ~~~ 0 < \alpha \leq 2, 
\end{equation}
and is obtained from the diffusion equation (or from the wave equation) by replacing the time derivative with a fractional derivative, in the Caputo sense.
\\
%\medskip
This explicitly reads in  the one-dimensional case with  ($D =1$) taking $u=u(x,t)$ and 
$\alpha \in (0,2]$
\begin{equation}
\frac{1}{\Gamma (n-\alpha)}\int_{0}^{t}{\frac{u^{(n)}(x,\tau)}{(t - \tau)^{\alpha + 1 - n}}\, d\tau} = 
\frac{\d^2}{ \d  x^2} u(x, t)\,,
\end{equation}
     where the derivatives $u^{(n)}$ are with respect time and  $n=1$ if $\alpha \in (0,1]$ and $n=2$ 
     if $\alpha \in (1,2]$.
     We recognize the for $\alpha \in (0,1)$ we have a fractional diffusion equation and when $\alpha\in (1,2)$
     we have a transition from the standard diffusion equation to the standard wave equation, that we refer to the fractional diffusion-wave equation.
      Of course, we find the standard diffusion and wave equations when 
     $\alpha=1$ and $\alpha=2$, respectively. 
%%%%%%%%%
\subsection{Cauchy and Signaling Problems}
The two basics problems for the Fractional PDE are the \emph{Cauchy} and \emph{Signaling} ones.
%\medskip
Denoting by $f(x)$ and $h(t)$ two given, well behaved functions:
\begin{equation}
\mathrm{Cauchy ~Problem:} 
\bigg \{
\begin{array}{rl}
u(x, 0^+) = f(x),& ~~ -\infty < x < +\infty \\
u(\pm \infty, t) =0, ~&~~~~~~~~~~~ t >0.\\
\end{array}
\end{equation}	
\begin{equation}
\mathrm{Signaling ~Problem:} 
\bigg \{
\begin{array}{rl}
u(x, 0^+) = 0, ~~~&~~~~~~ 0 < x < +\infty \\
u(0^+, t) =h(t),& u(+\infty , t) = 0, ~ t > 0 .\\
\end{array}
\end{equation}
%\medskip
When $1 < \alpha \leq 2 $ an extra initial condition is needed for both problems: $u_t(x,0^+) = g(x)$.
\\
The solutions turn out to be expressed by  proper convolutions between the source function and a (two) characteristic function(s), the so called \emph{Green functions} or \emph{fundamental solutions} $\mathcal{G}(x,t)$ of the problem.\\
The function $\mathcal{G}_C(x,t)$ for the Cauchy problem represents the solution for $f(x) = \delta (x)$, and the function $\mathcal{G}_S(x,t)$ for the Signaling problem represents the solution for $h(t) = \delta _+ (t)$.\\
The Green functions are connected by the {\it Reciprocity Relation}:
\begin{equation}
2\nu x \mathcal{G}_C(x,t; \nu) = t \mathcal{G}_S(x,t; \nu) =  \nu z M_{\nu}(z),
\end{equation}
with $\nu = {\alpha}/{2}$ and $z = {x}/{t^{\nu}}$ being the {\it similarity variable}.
Here  $M_\nu(z)$  is a function of the Wright type, entire  for $0\le \nu <1$, referred to as
 {\it $M$-Wright function} briefly  discussed below.
For more details see e.g the survey paper by Mainardi, Luchko and Pagnini \cite{Mainardi LUMAPA01}
based on the previous papers by Mainardi
 \cite{Mainardi AML96,Mainardi CHAOS96,Mainardi CISM97}. 
%%%%%%%%    
\subsection{The Fundamental Solution of the Cauchy problem}
According to Eq. (13) the fundamental solution of the Cauchy   problem is thus provided by
\begin{equation}
\G_C(x,t;\nu) = \frac{1}{2 t^\nu}\, M_\nu\left( \frac{x}{t^\nu}\right)\,,
\end{equation}
where the function $M_{\nu}(z)$ reads in its series and integral representations as
\begin{equation}
M_{\nu}(z) := \sum_{n=0}^{\infty}{\frac{(-z)^n}{n!\mathit{\Gamma}[-\nu n + (1-\nu)]}} = 
\frac{1}{2\pi i}{\int_{Ha}{e^{\sigma - z\sigma ^{\nu}}\frac{d\sigma }{\sigma ^{1-\nu }}}},
\end{equation} 
where $Ha$ denotes the Hankel path.
%  a contour that begins at $t = -\infty - ib$  $(b > 0$), encircles the branch cut that lies along the
% negative  real axis, and ends up at$ t = -infty 1+ ia$ $(a > 0$).
$M_{\nu}(z)$
 results to be a particular case of the Wright function  $W_{\lambda , \mu}(z)$:
\begin{equation}
W_{\lambda , \mu}(z) := \sum_{n=0}^{\infty}{\frac{z^n}{n! \mathit{\Gamma}(\lambda n + \mu)}}, ~~~\lambda > -1, ~ \mu \in \mathbb{C}, ~ z \in \mathbb{C}
\end{equation}
when $\lambda =-\nu$ and $\mu =1-\nu$.
%\\
We note that for convenience and for historical reasons the Wright functions may be classified in two kinds: the first kind with $\lambda\ge 0$ and the second kind with $-1<\lambda <0$  so following the Appendix F 
of Mainardi's book \cite{Mainardi BOOK10}, where the reader can find also some historical notes.  
%\\
\begin{figure}%
	\centering
%	\subfloat
{{\includegraphics[width=7.0cm]{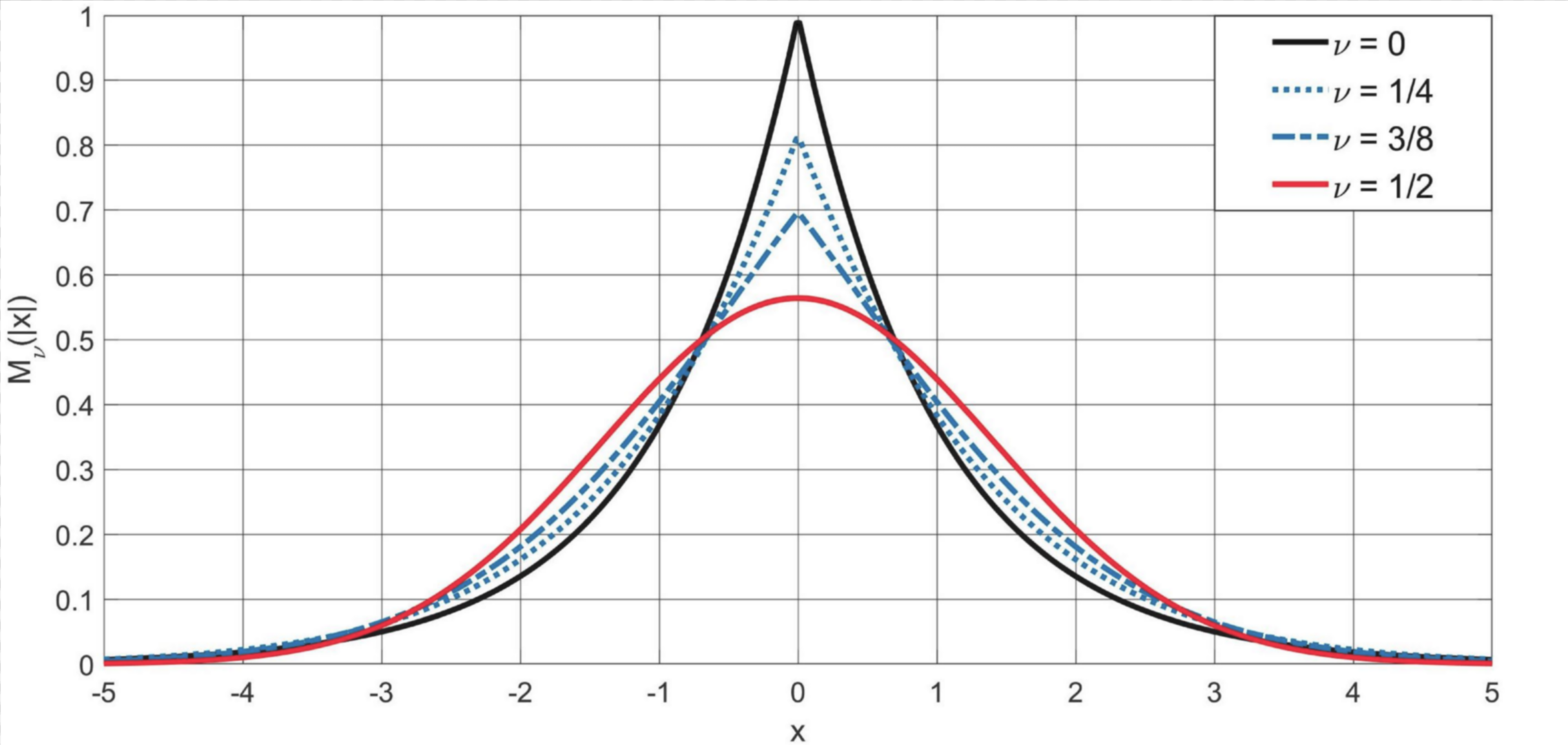}}}%
	\quad
%	\subfloat
{{\includegraphics[width=7.0cm]{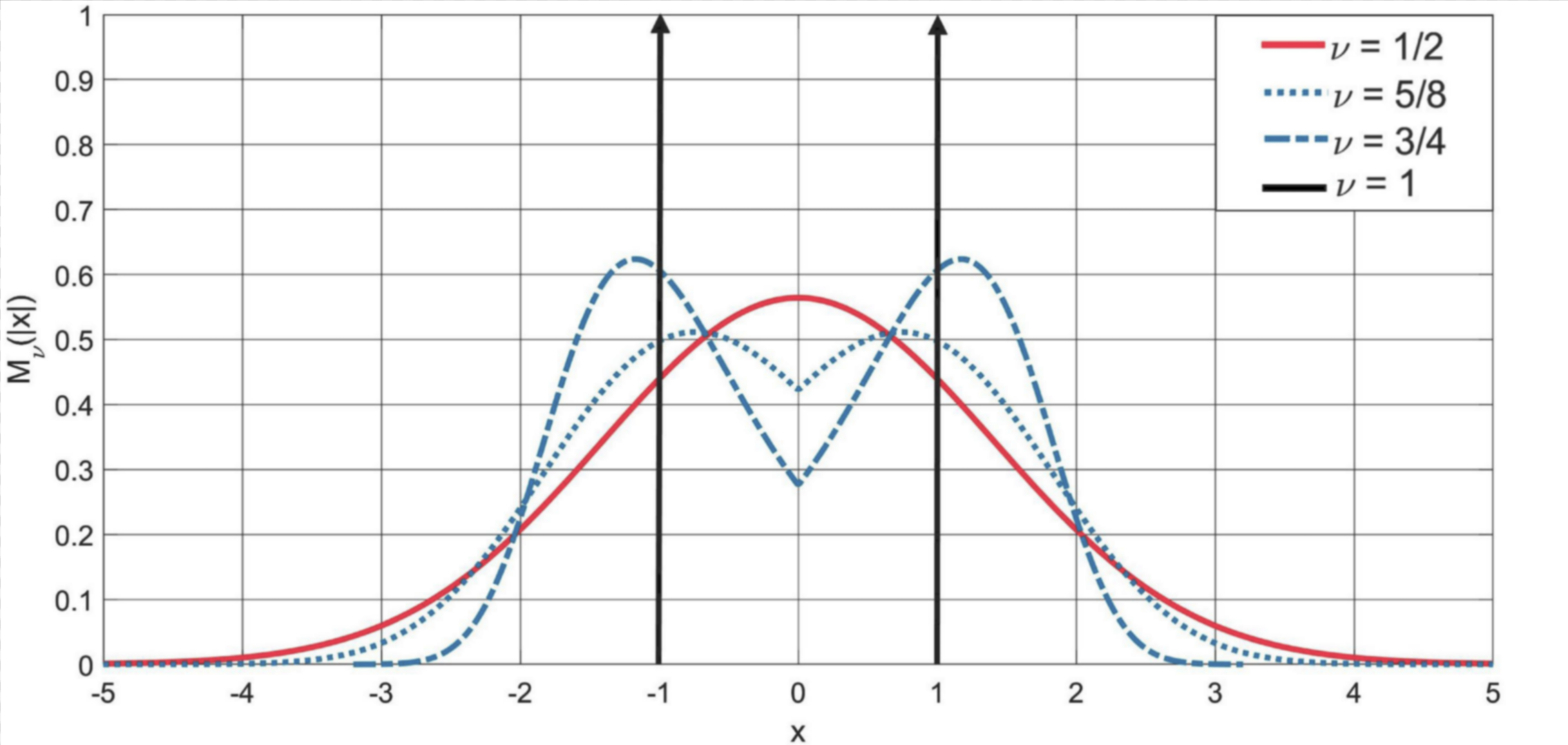}}}%
	\caption{Plots of the functions $M_{\nu}(|x|)$  for $|x| \le 5$ 
	 at $ t=1$;
	 LEFT: for $\nu = 0, 1/4,  3/8, 1/2$.
	 RIGHT: for  $\nu = 1/2, 5/8, 3/4, 1$}
	\label{fig:mainardifunction}%
\end{figure}
%\\
For $\nu = 1/2$, we find back the diffusion equation and indeed the $M_{\nu}(|x|)$ function becomes the Gaussian function known as fundamental solution of the diffusion equation for  the Cauchy Problem) whereas , for $\nu \rightarrow 1$, we find back the wave equation and indeed the $M_{\nu}(|x|)$ function tends to two Dirac delta functions ad fundamental solutions of the Cauchy Problem, 
centered in $x = \pm 1$.
\vsp
It should be noted that in the 1993  book by Pr{\"u}ss \cite{Pruss BOOK93}  we find a figure 
quite similar to our Fig. 5 reporting the $M$-Wright function in
for $\nu \in [1/2, 1)$. 
% namely the Wright function of the second kind. 
It was  derived from inverting the Fourier transform expressed in terms of 
the Mittag-Leffler function following the approach 
by Fujita \cite{Fujita 90} for the fundamental  solution of the Cauchy problem for 
the diffusion-wave equation, fractional in time. Both approaches were carried out  without any relation to the Wright function,
 presumably unknown to Fujita and Pr{\"u}ss.
%\\
However, our plot must be considered independent from  that of Pr{\"u}ss because 
Mainardi  
used the Laplace transform in his former paper presented at WASCOM, Bologna, 
October  1993
\cite{Mainardi WASCOM93}
(and published later in a number of papers and in his 2010 book \cite{Mainardi BOOK10})  
so he  was  aware of the book by Pr{\"u}ss only later. 
%%%%%%%%%%%%
%\vsp
%%%%
Furthermore, we invite  readers to look at the simulation of the fundamental solution
$M_\nu(x,t)$ at $t=1$ for $x \in [-5,+5]$   at varying $\nu$ available in  YOUTUBE:  
%% https://www.youtube.com/watch?v=uf_4aB1COPg
{\texttt{https://www.youtube.com/watch?v=uf{\_}4aB1COPg}}
carried out by Consiglio with Matlab.
%%%%%%%%%%
\subsection{Plots of Solutions for the Cauchy Problem}
We start this graphical sub-section by reporting the plots of the evolution in time of fundamental solutions  (that is $f(x) = \delta(x)$ and $g(x)=0$) taking $\nu=0.65$, $\nu=0.75$ and $\nu=0.85$.
\\
The general solutions for the Cauchy problem  are  obtained through suitable  convolution in space as shown hereafter:
\begin{equation}
u(x,t;\nu) = \int_{-\infty}^{+\infty}{\mathcal{G}_C(\xi,t; \nu)f(x-\xi)d\xi}, ~~~~ 0 < \nu \leq 1/2,
\end{equation}
\begin{equation}
\begin{split}
u(x,t;\nu) = \int_{-\infty}^{+\infty}{\biggl[ \mathcal{G}^{(1)}_C(\xi,t; \nu)f(x-\xi)   + \mathcal{G}^{(2)}_C(\xi,t; \nu)g(x-\xi) \biggl]d\xi  },\\ ~~~~~~~~~~~~~~~~~~~~~~~~~~~~~~~~~~~~~~~1/2 < \nu \leq 1,
\end{split}
\end{equation}
where $\mathcal{G}^{(2)}_C$ is the primitive in time of $\mathcal{G}^{(1)}_C$.
%For $\nu = 0.65$ and $f(x) = \delta (x)$:
\begin{figure}[h!]
	\centering{{
			\includegraphics[width=5.5cm]{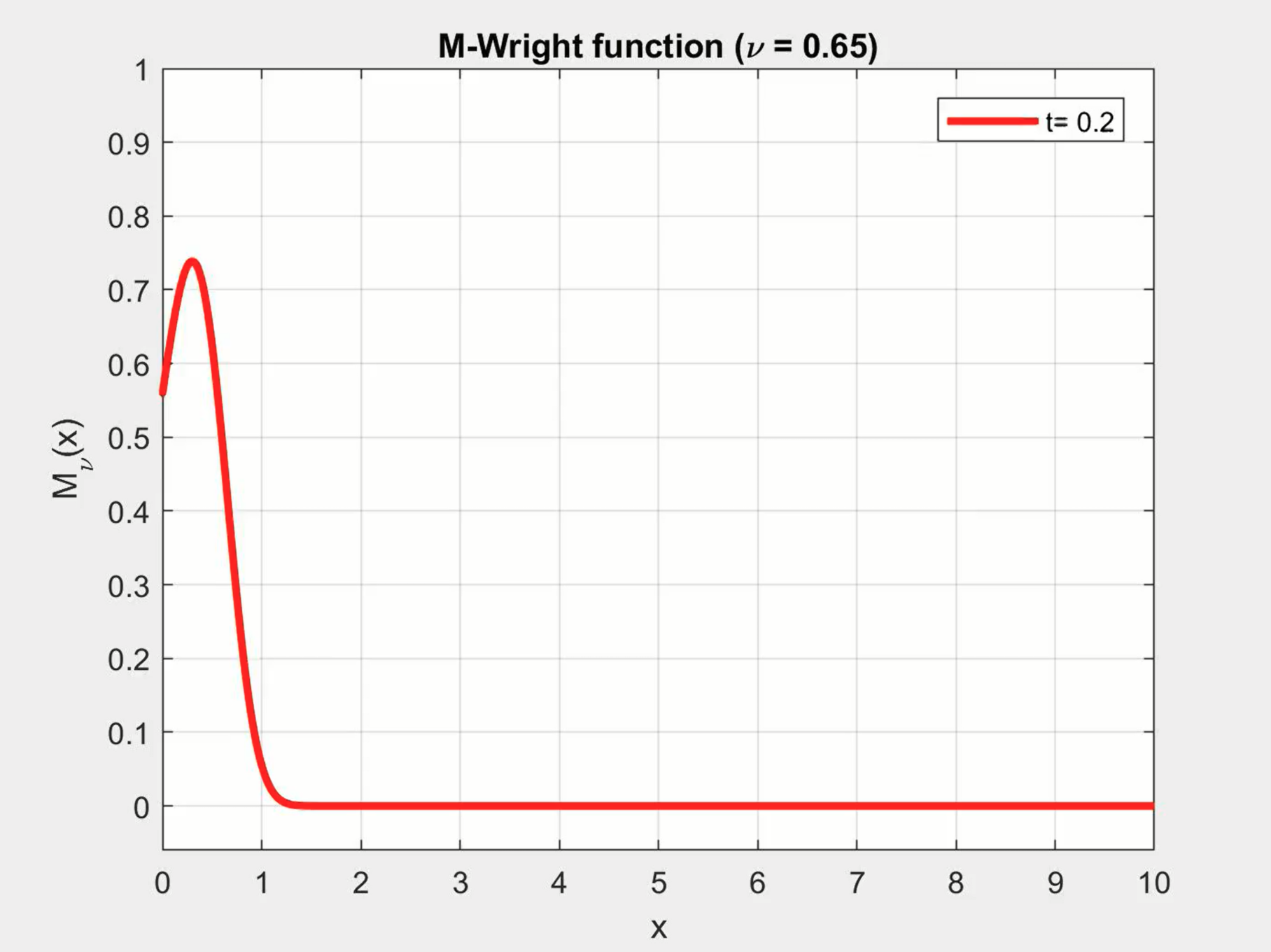}
		}
		{
			\includegraphics[width=5.5cm]{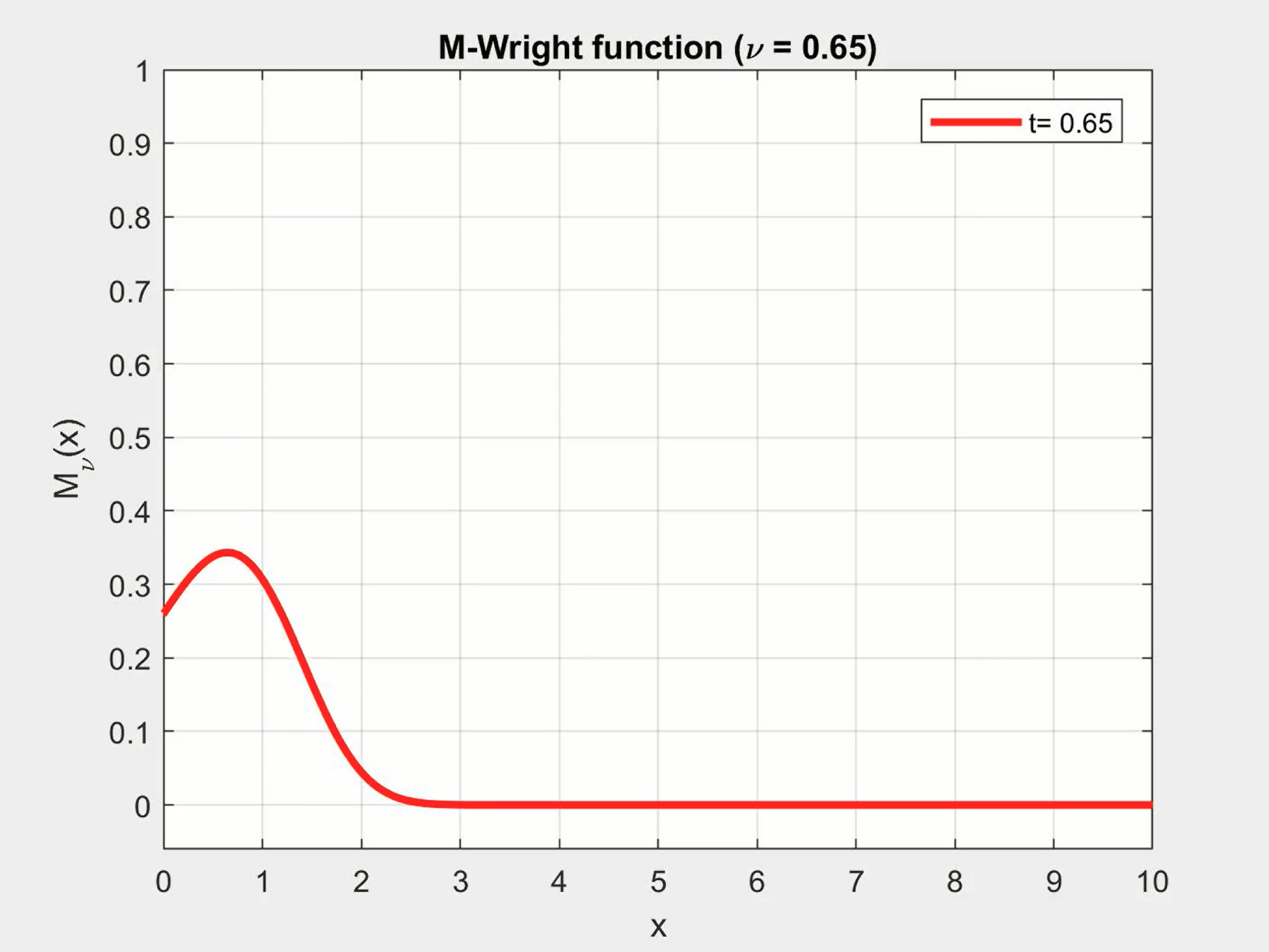}
		}
		\\
		{
			\includegraphics[width=5.5cm]{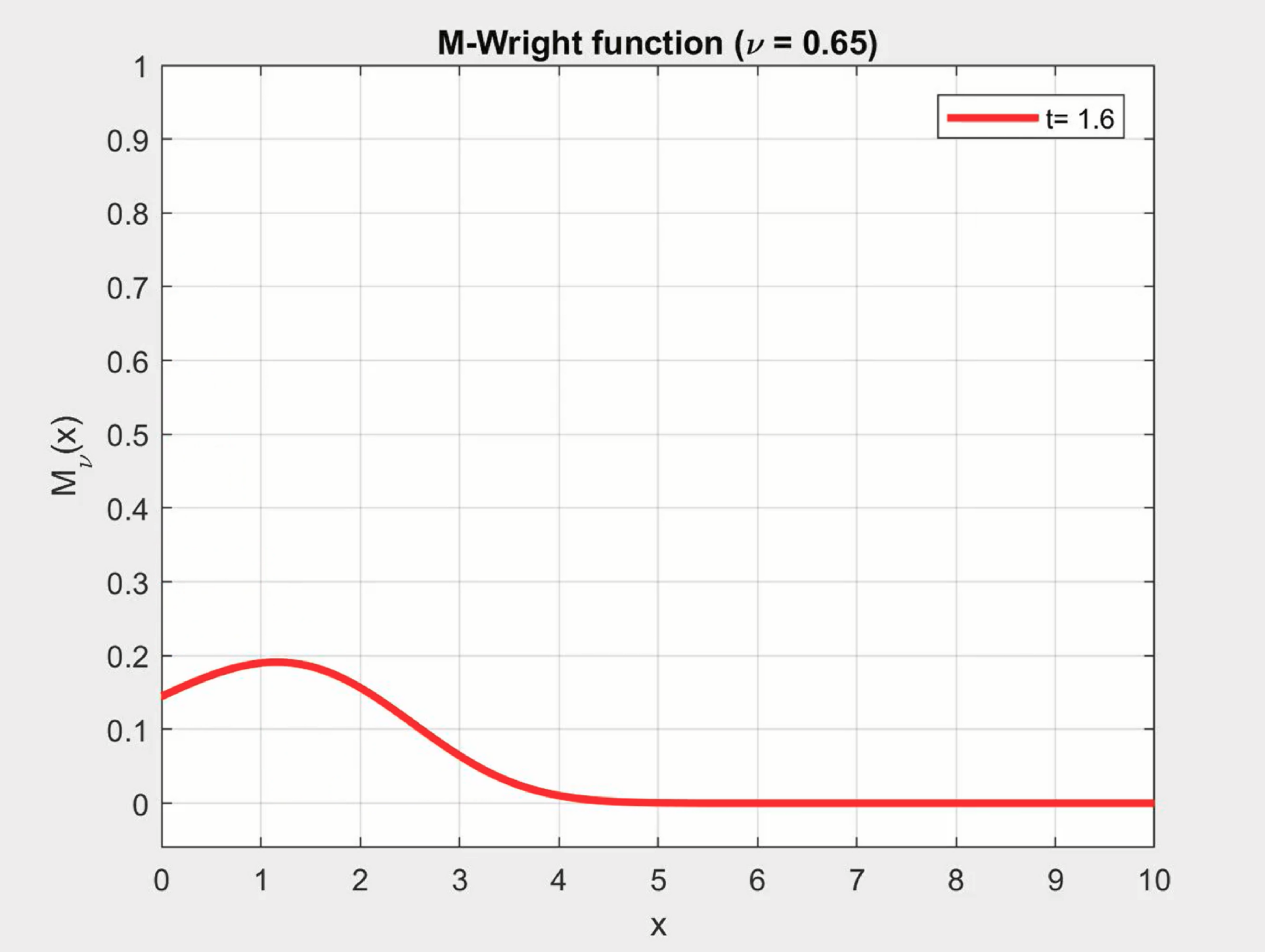}
		}
		{
			\includegraphics[width=5.5cm]{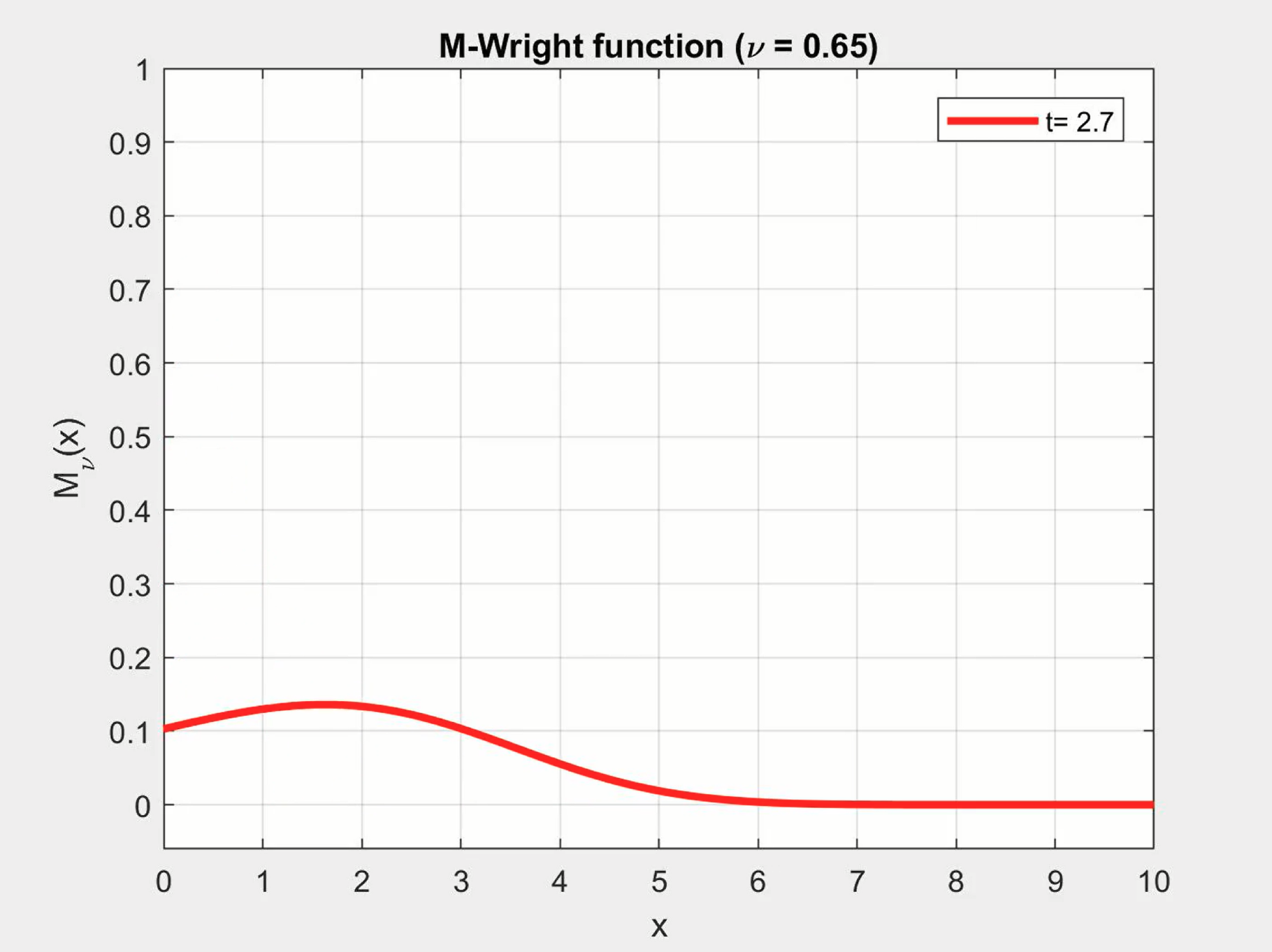}
	}}
	\caption{Time evolution of the fundamental solution for $ \nu =0.65 $ in the Cauchy problem}
\end{figure}
%%%
\vskip 3truecm
%For $\nu = 0.75$ and $f(x) = \delta (x)$:
\begin{figure}[h!]
	\centering{{
			\includegraphics[width=5.5cm]{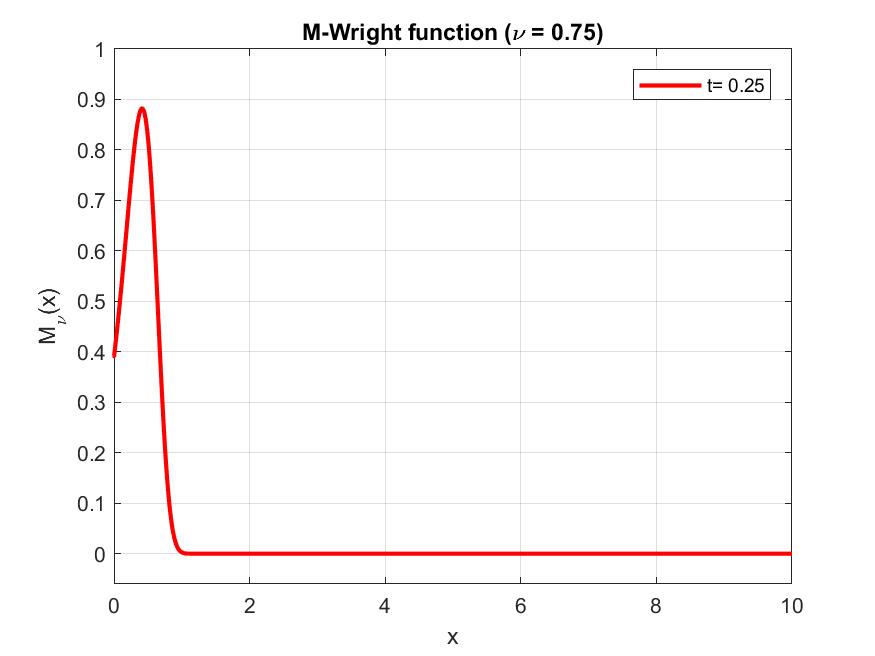}
		}
		{
			\includegraphics[width=5.5cm]{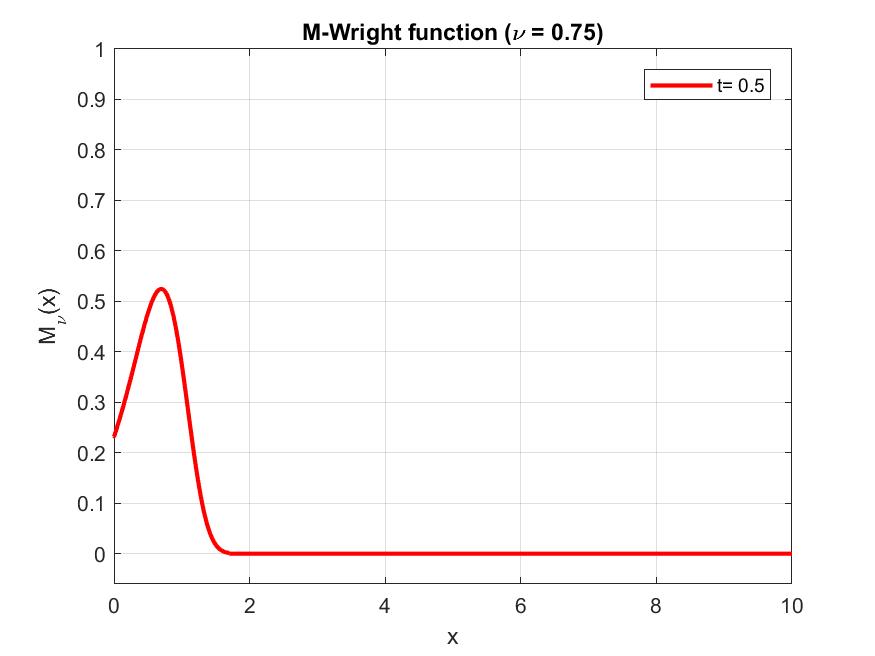}
		}
		\\
		{
			\includegraphics[width=5.5cm]{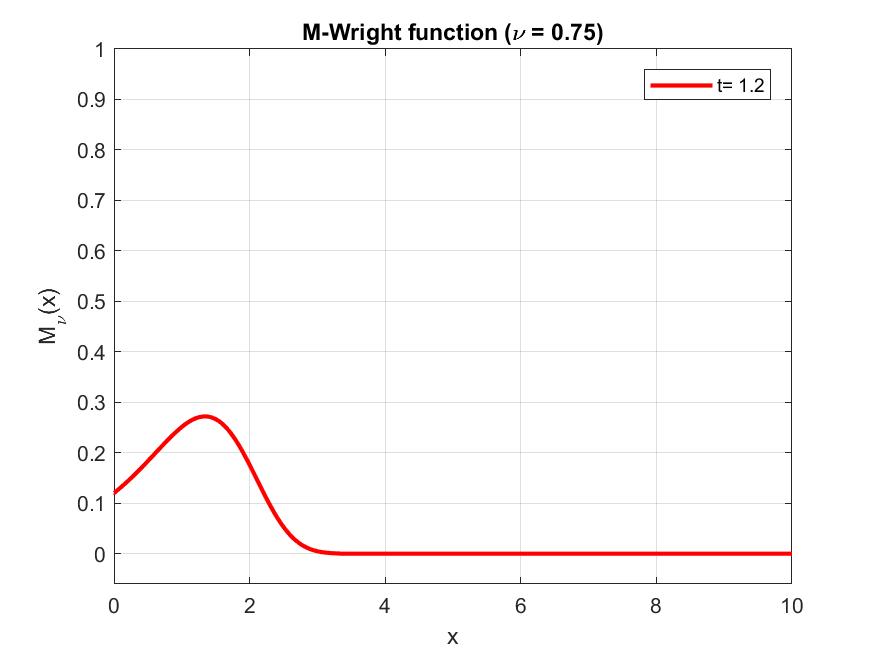}
		}
		{
			\includegraphics[width=5.5cm]{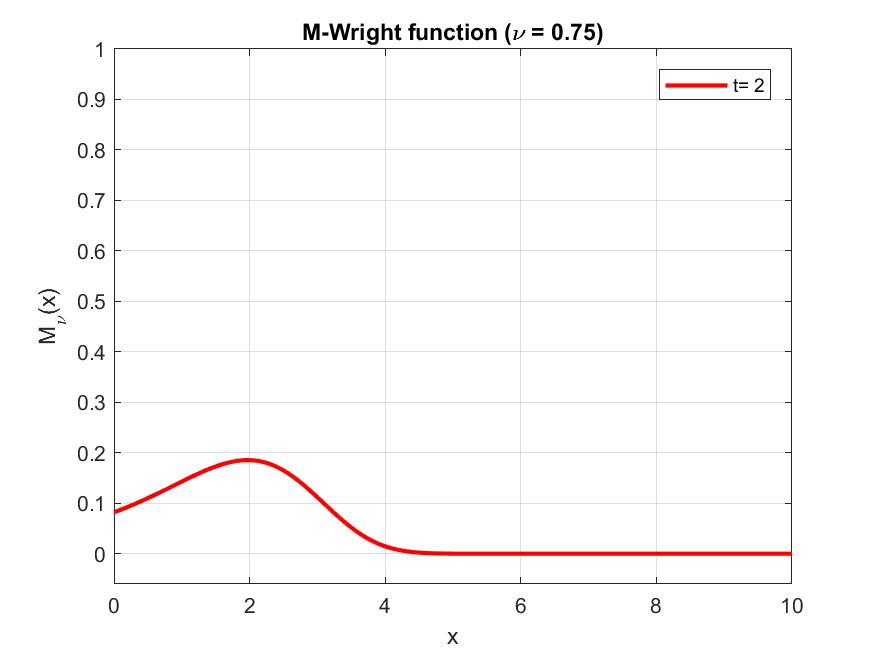}
	}}
	\caption{Time evolution of the fundamental solution for $\nu=0.75$ in the Cauchy problem}
\end{figure}
\vskip 2truecm
%%For $\nu = 0.85$ and $f(x) = \delta (x)$:
\begin{figure}[h!]
	\centering{{
			\includegraphics[width=5.5cm]{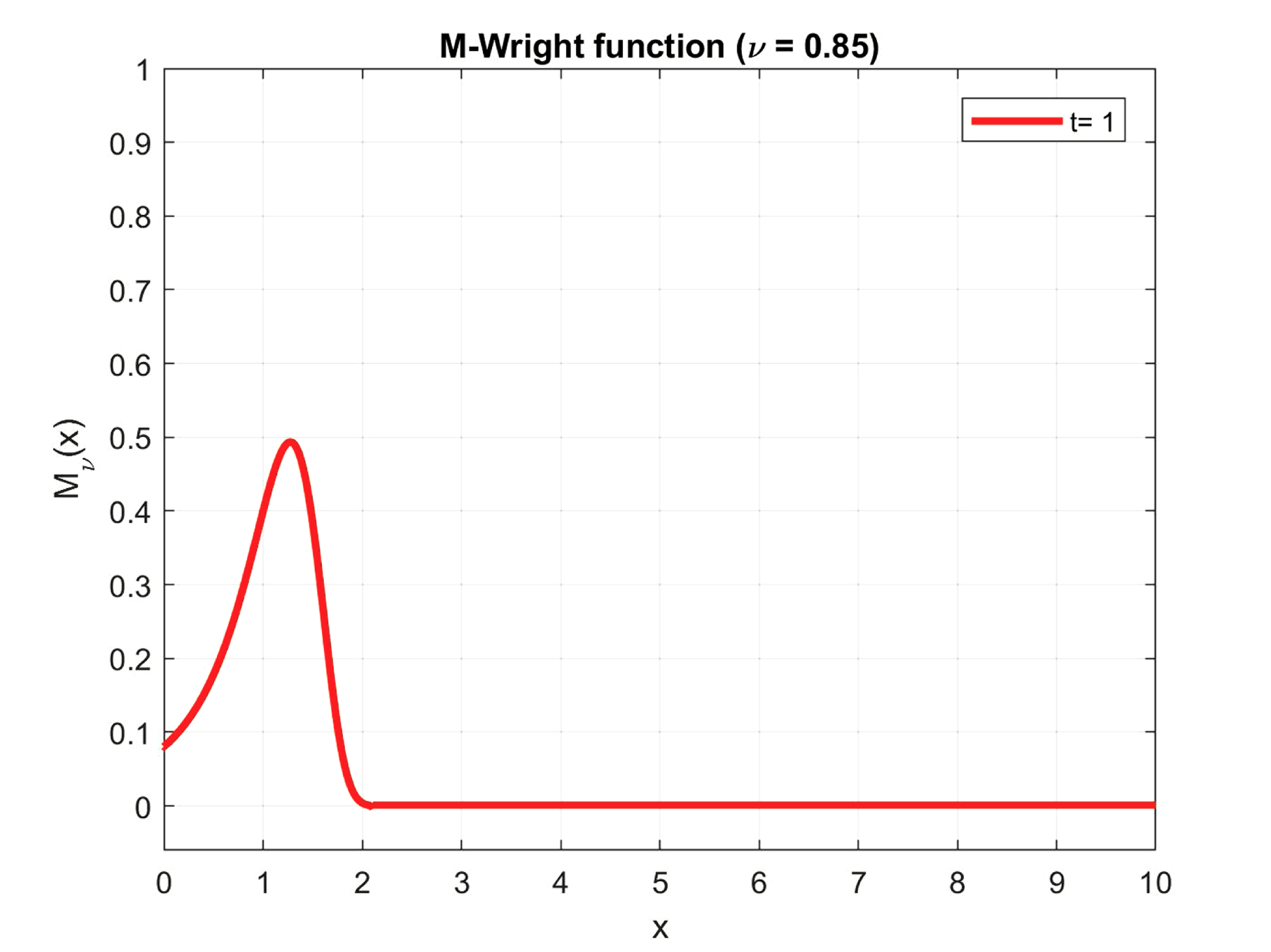}
		}
		{
			\includegraphics[width=5.5cm]{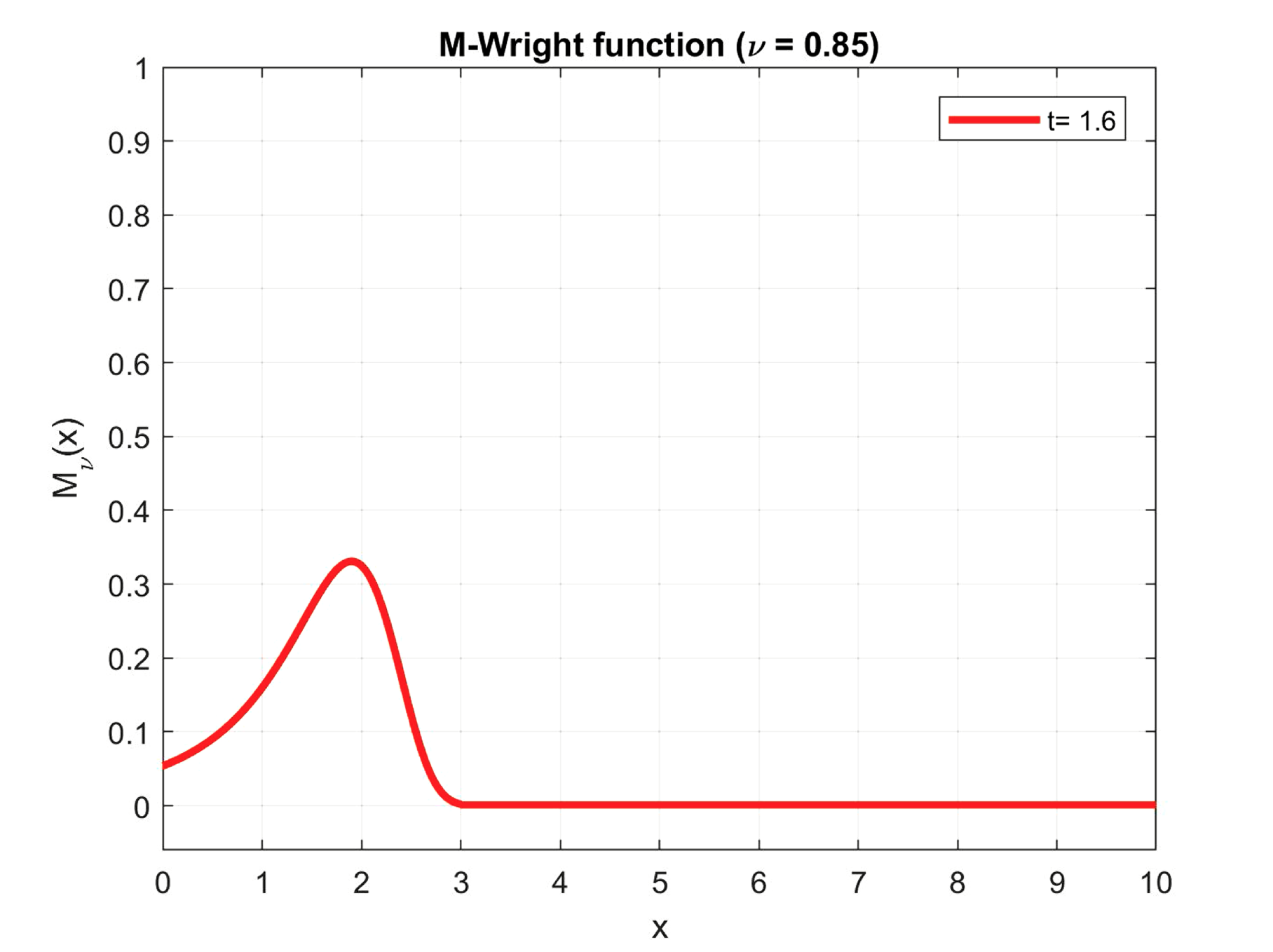}
		}
		\\
		{
			\includegraphics[width=5.5cm]{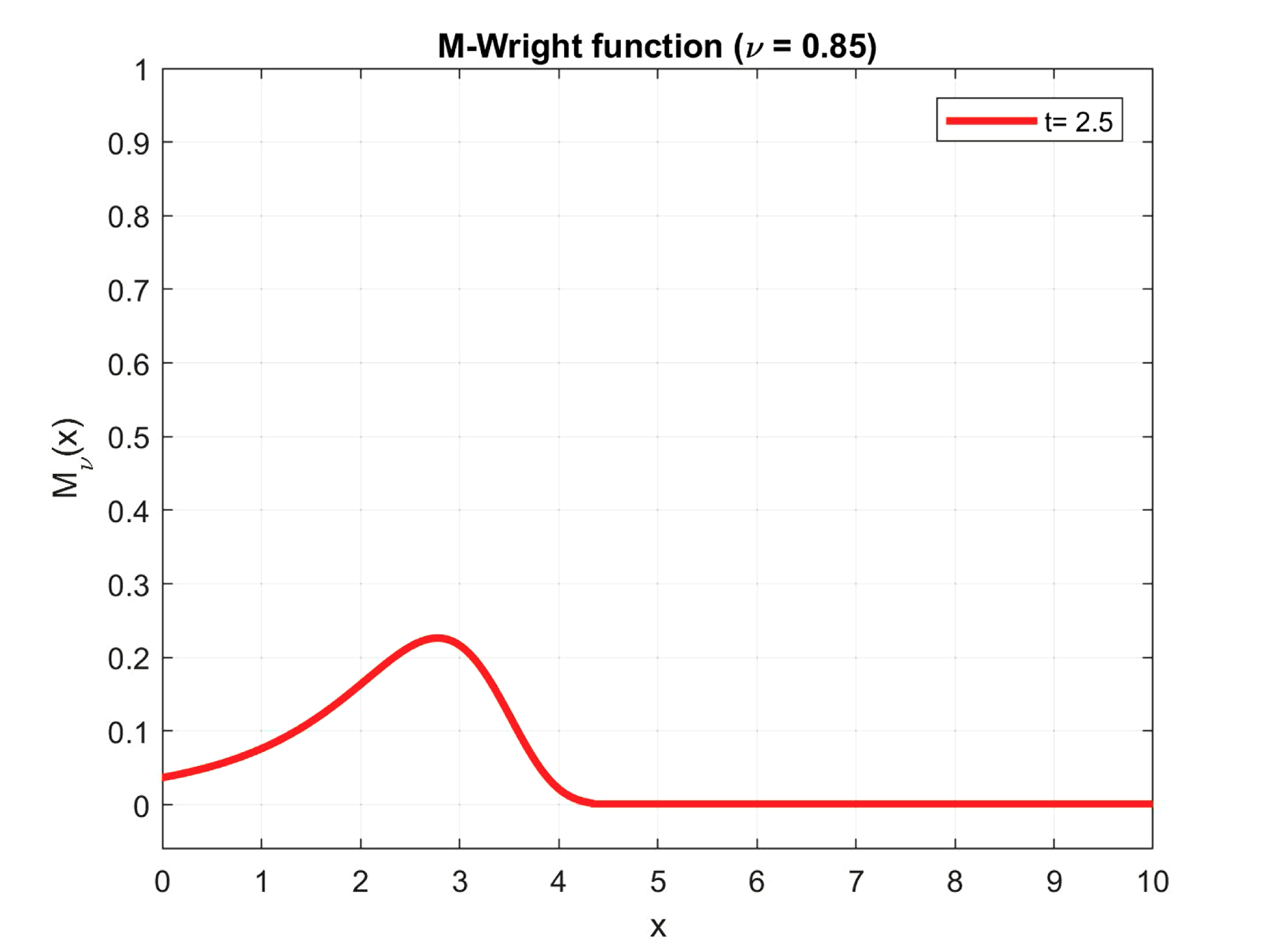}
		}
		{
			\includegraphics[width=5.5cm]{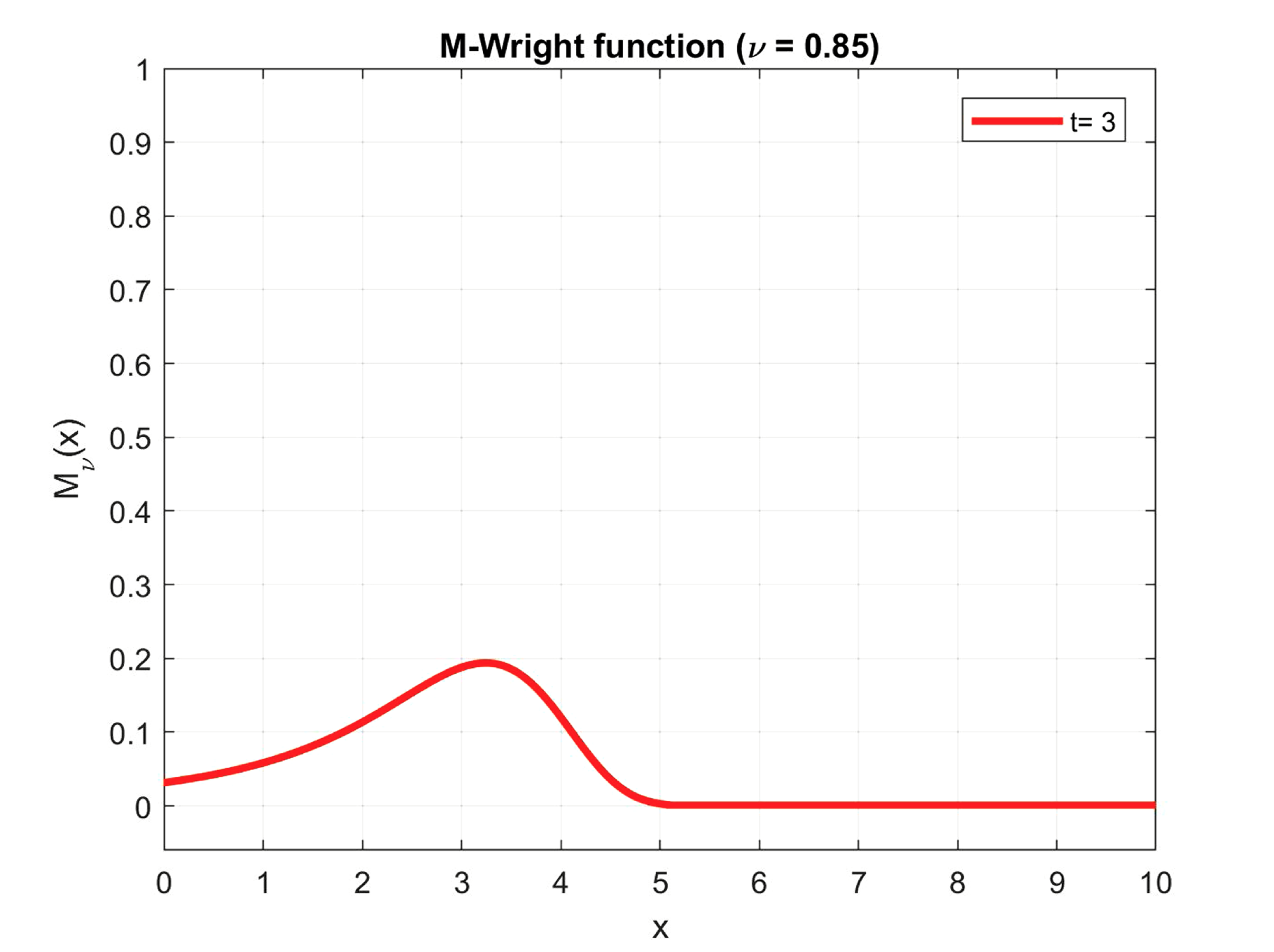}
	}}
	\caption{Time evolution of the fundamental solution for $\nu =0.85$ in the Cauchy problem}
\end{figure}
%\newpage 
%%%%%%%%%%

%\\
%\medskip
%\newpage
\vsp
We now consider
  for  the Time-Fractional Diffusion Wave equation
the evolution of an initial centered box-signal $f(x) =1$ for $-1 \le  x \le +1 $ and zero otherwise, 
assuming  $g(x) = 0$.
%$f(x) = H(1-|x|)$ where $H(x)$ is the Heaviside step function.
  Based on Eqs.  (17)- (18)
we simulate the cases $\nu=0.50$ (standard diffusion),  $\nu=0.75$ and $\nu=1$ (standard wave), as shown in the next figures where the initial box function is denoted by dotted lines.
%\\
\vskip-0.5truecm
\begin{figure}[h!]%
	\centering
%	\subfloat
{{\includegraphics[width=5.5cm]{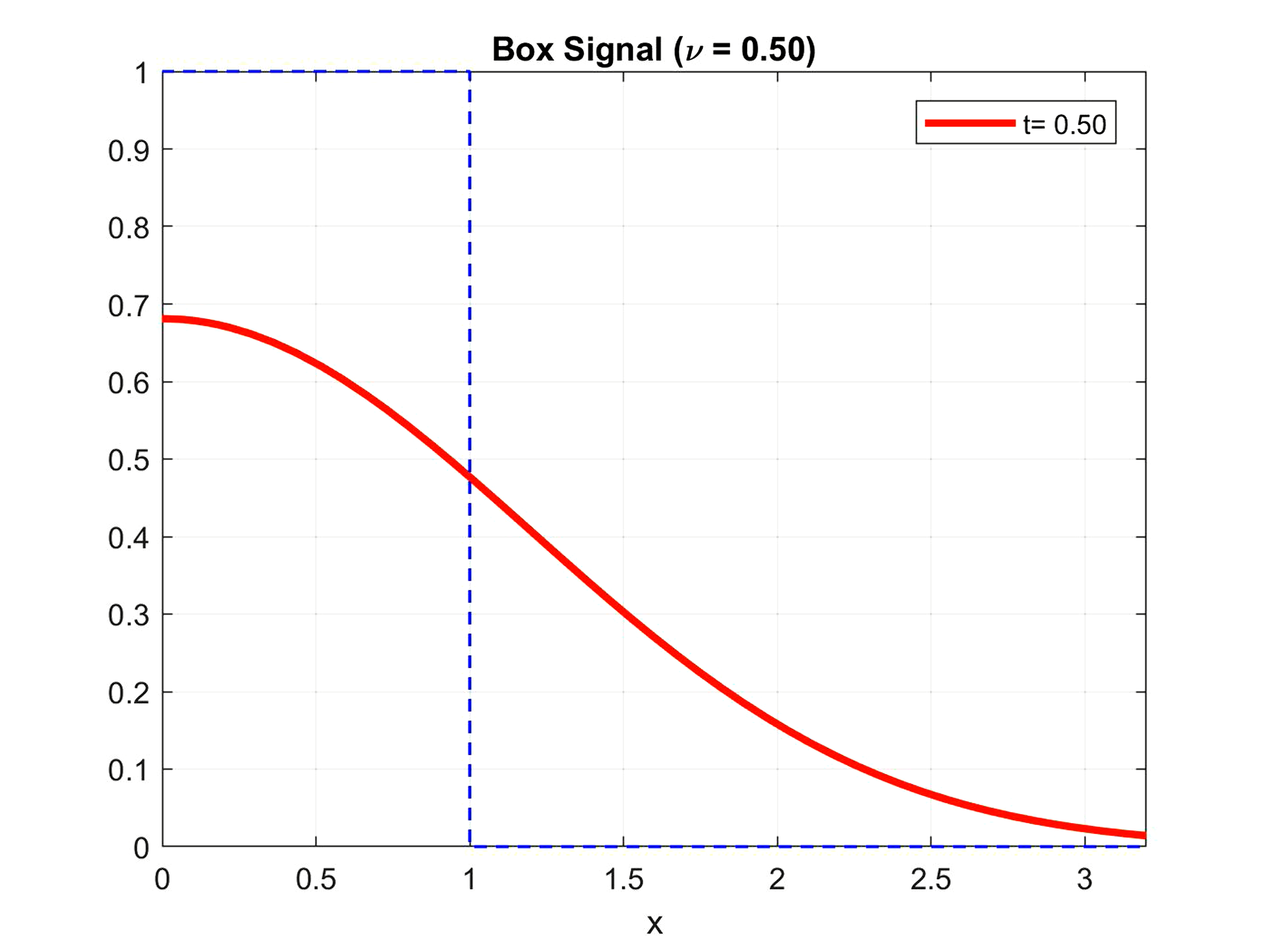} }}%
	\qquad
%	\subfloat
{{\includegraphics[width=5.5cm]{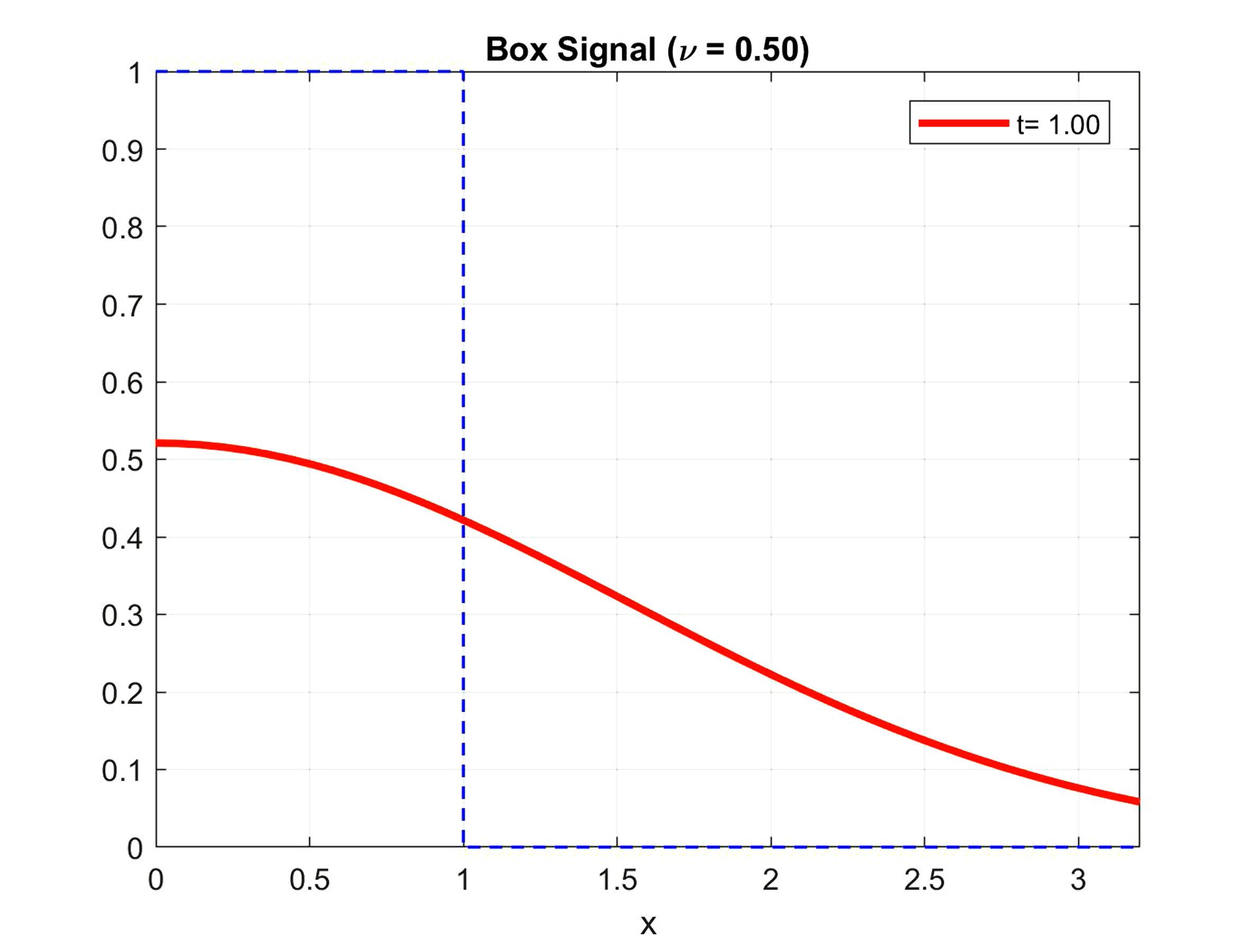} }}%
	\caption{Time evolution of an initial box-signal for $\nu=0.50$, seen at t = 0.50 (left) and t = 1.00 (right) in $0 \leq x \leq 3.5$. The problem is symmetric on the negative axis. }
	\label{fig:boxgauss}%
\end{figure}
\vskip-0.5truecm
%For $\nu = 0.75$ and $f(x) = H(1-|x|)$:
\begin{figure}[h!]%
	\centering
%	\subfloat
{{\includegraphics[width=5.5cm]{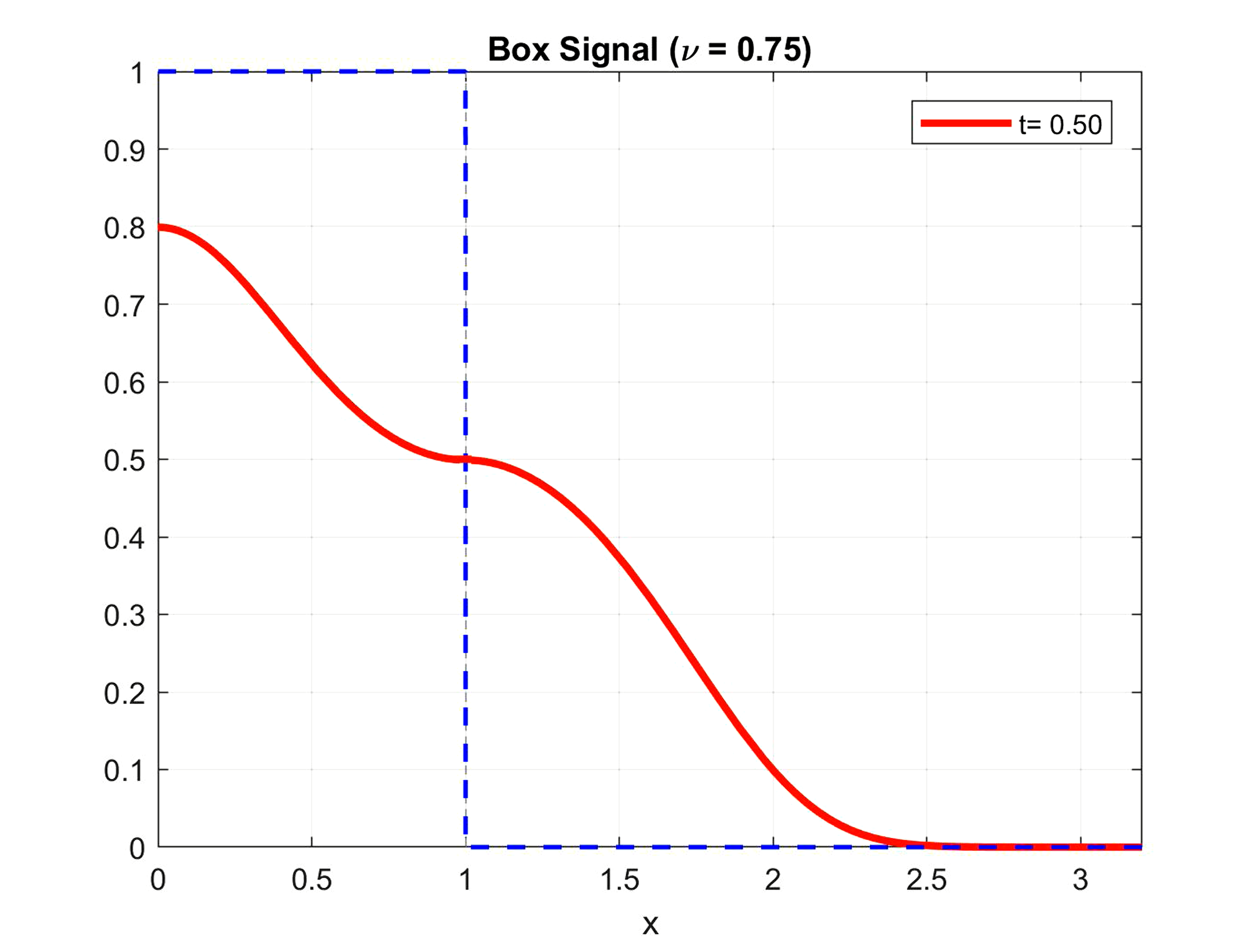} }}%
	\qquad
%	\subfloat
{{\includegraphics[width=5.5cm]{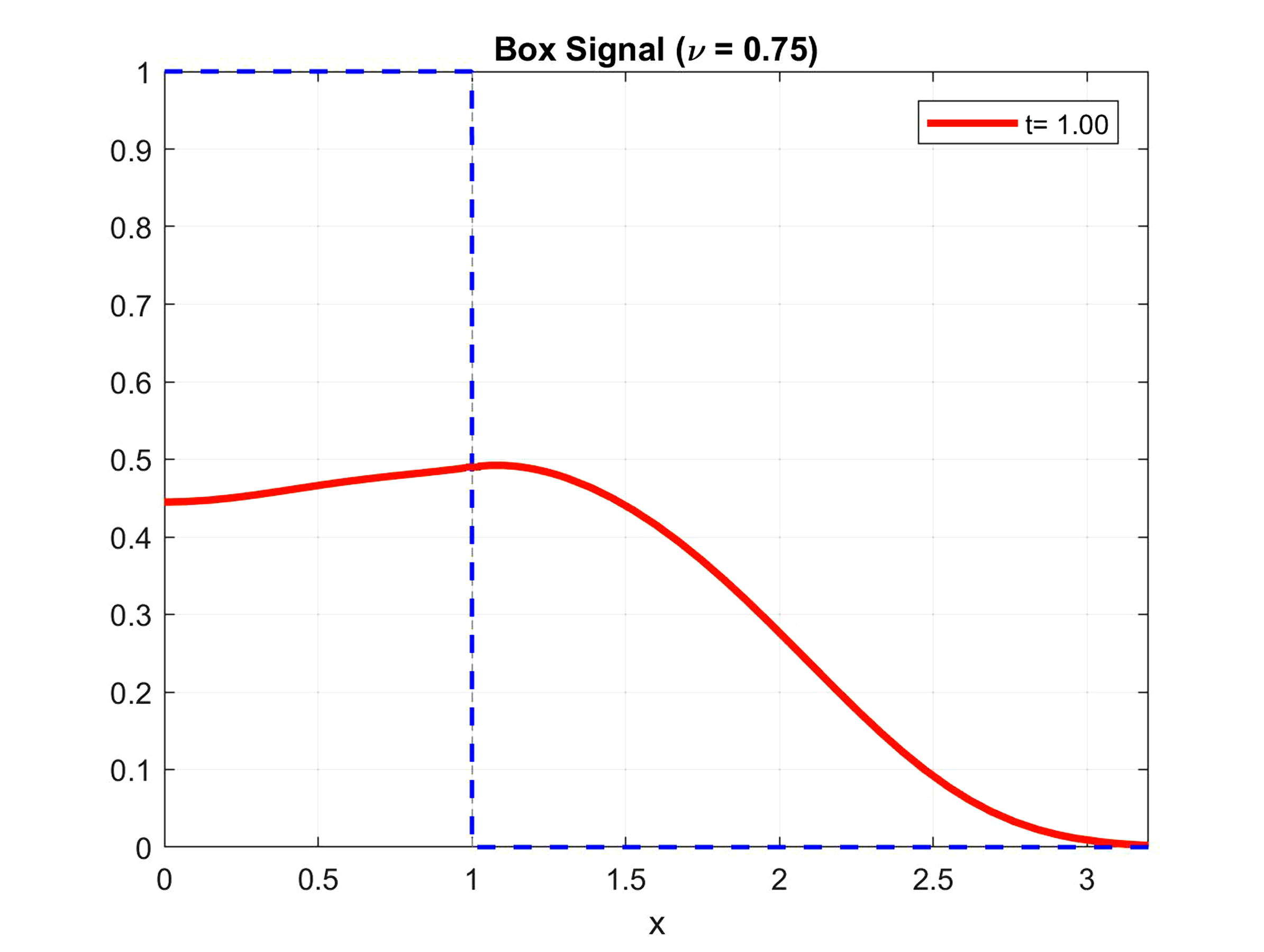} }}%
	\caption{Time evolution of an initial box-signal for $\nu=0.75$, seen at t = 0.50 (left) and t = 1.00 (right) in $0 \leq x \leq 3.5$. The problem is symmetric on the negative axis. }
	\label{fig:box75}%
\end{figure}\vskip-0.5truecm
%For $\nu = 1$ and $f(x) = H(1-|x|)$:
\begin{figure}[h!]%
	\centering
%	\subfloat
{{\includegraphics[width=5.5cm]{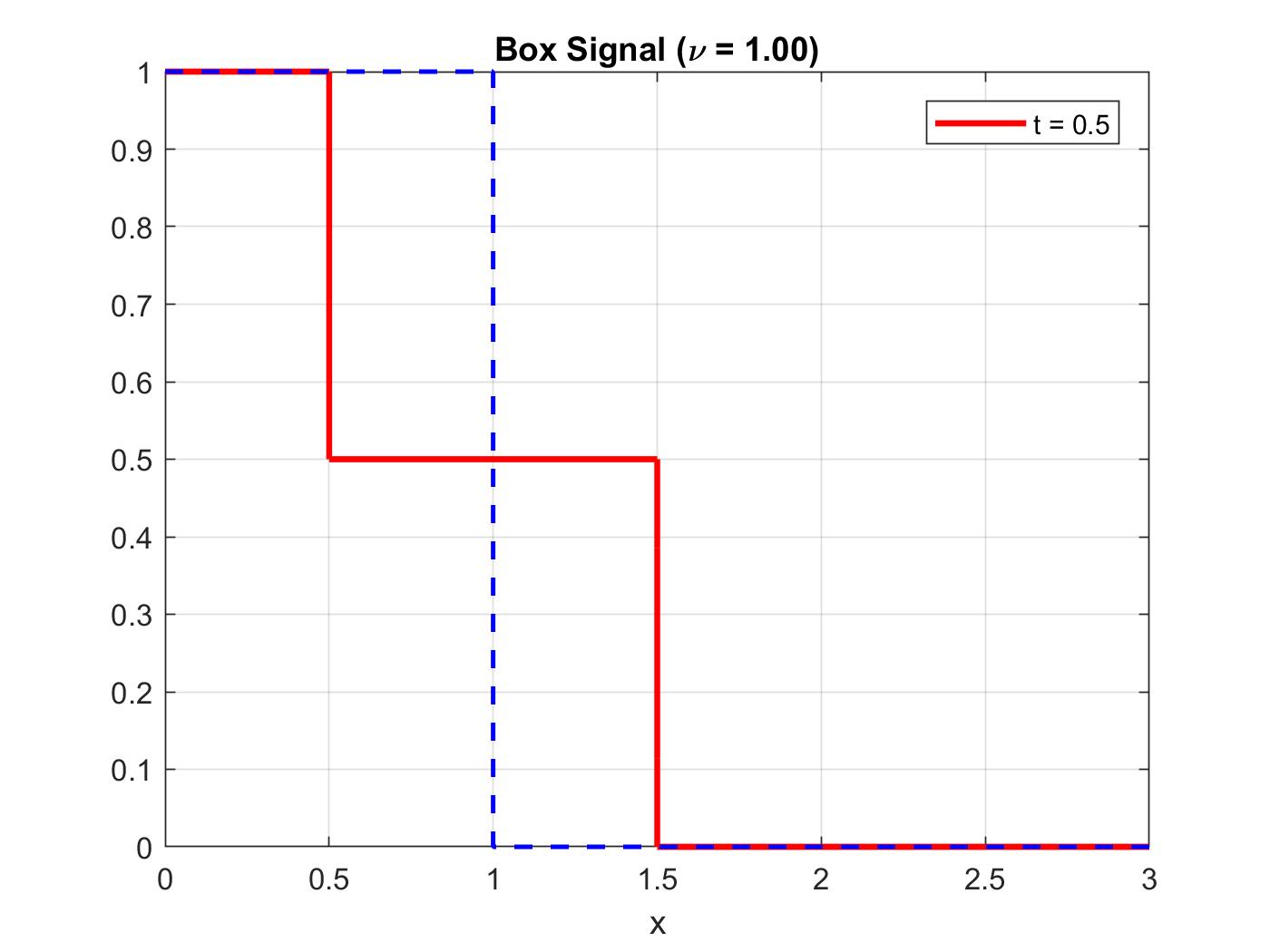} }}%
	\qquad
%	\subfloat
{{\includegraphics[width=5.5cm]{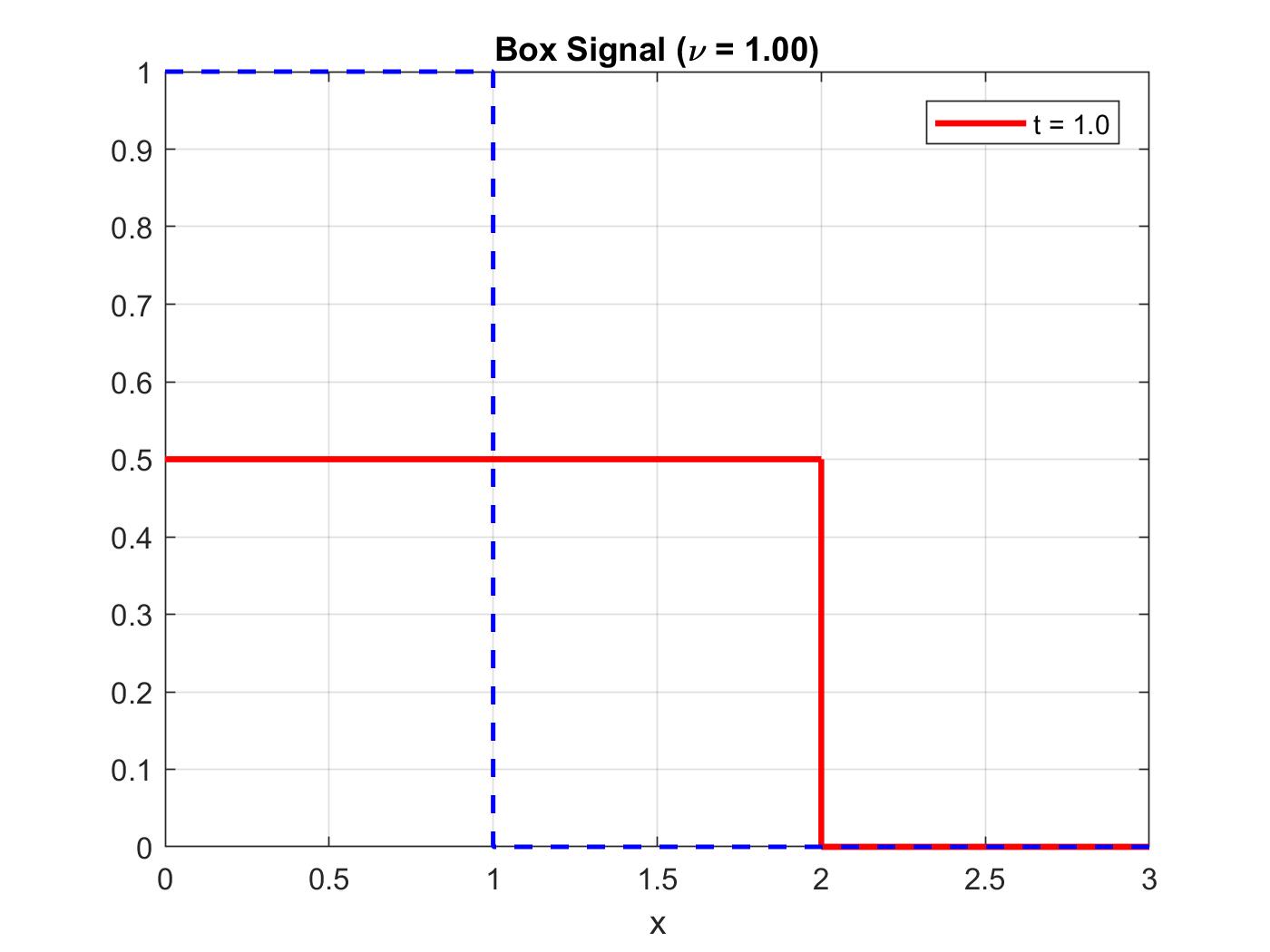} }}%
	\caption{Time evolution of an initial box-signal for $\nu=1$, seen at t = 0.50 (left) and t = 1.00 (right) in $0 \leq x \leq 3$. The problem is symmetric on the negative axis. }
	\label{fig:boxwave}%
\end{figure}
%%%%%%%%%%%
\newpage
 \section{Concluding remarks}
%%%%%%%%%%%%    
    We have considered the  so-called time  fractional diffusion wave equation with particular attention when this equation is interpolating the processes of diffusion and the wave propagation.
    In these cases we speak about fractional diffusive waves using a term incorporating both diffusion and wave phenomena.
    \\
    We have analyzed  and simulated  both the situations in which the input  function is a Dirac delta generalized function and a box function, restricting ourselves to the Cauchy problem.
 In the first case we get  the fundamental solutions
(or Green function) of the problem whereas  in the latter case the solutions are obtained by a 
space convolution of the Green function with the input  function.
\\
In the next future we plan to consider the signaling problem in order to complete the topic of fractional diffusive waves.
%%%%%%%%%
\section*{Acknowledgments}
The work of FM  has been carried out in the framework of the activities of the National Group of Mathematical Physics (GNFM, INdAM).
\\
The results contained in the present paper have been
partially presented at the XX International Conference {\it Waves and Stability in Continuous Media} (WASCOM 2019)
held in Maiori (Sa),  Italy, June 10-14 (2019).
  %%%%%%%%%
    
%%%%%%%%%%%
\end{document}